\documentclass[12pt]{article}
\usepackage{amsfonts}
\usepackage{amsmath}
\usepackage{array}
\usepackage{listings}
\usepackage{epsfig}
\usepackage{xcolor}
\usepackage{geometry}
\title{A Universal $8$-Vertex Origami Torus Triangulation}

\author{Peter G. Doyle, Fabian Lander%
\thanks{Supported by the International Max Planck Research School Mathematics in the Sciences}, \\Richard Evan Schwartz%
\thanks{Supported by N.S.F. Research Grant DMS-2505281}, Steve Trettel}

\newtheorem{theorem}{Theorem}[section]
\newtheorem{proposition}[theorem]{Proposition}
\newtheorem{lemma}[theorem]{Lemma}

\newtheorem{corollary}[theorem]{Corollary}

\newcounter{fig}[section]
\renewcommand{\thefig}{\thesection.\arabic{fig}}
\newcommand{\figcaption}[2]{\par\vspace{8pt}\refstepcounter{fig}\label{#1}{\bf Figure \thefig:\/} #2}

\newtheorem{remark}[theorem]{Remark}

\def\startproof{{\bf {\medskip}{\noindent}Proof: }}

\def\endproof{$\spadesuit$  \newline}

\def\C{\mbox{\boldmath{$C$}}}%
\def\H{\mbox{\boldmath{$H$}}}%
\def\R{\mbox{\boldmath{$R$}}}%
\def\Z{\mbox{\boldmath{$Z$}}}%

\begin{document}

\maketitle

\begin{abstract}
  An origami torus is a piecewise affine isometric  embedding of a flat
  torus into $\R^3$.   We prove that every flat torus is realized
  as an $8$-vertex origami torus based on the triangulation
  with degree sequence $4,6,6,6,6,6,7,7$.   This provides an
  $8$-vertex universal triangulation for origami tori.
    \end{abstract}

\section{Introduction}
\label{sec:intro}

\subsection{The Main Result}

A {\it flat torus\/} is a quotient of the
form $T=\R^2/\Lambda$, where $\Lambda$ is a
lattice of translations of $\R^2$.  The lattice
$\Lambda$ is called the {\it lattice\/} of $T$.
A {\it triangulation\/} ${\cal T\/}$ of a flat torus $T$ is
a decomposition of $T$ into finitely many triangles,
such that every two triangles are either disjoint,
or share one vertex, or share one edge. 

Two triangulations ${\cal T\/}_1$ and ${\cal T\/}_2$, respectively
of flat tori $T_1$ and $T_2$, are {\it combinatorially
  equivalent\/} if there is a
homeomorphism $h: T_1 \to T_2$
which affinely maps each triangle of ${\cal T\/}_1$ to a triangle
of ${\cal T\/}_2$.    The equivalence class of ${\cal T\/}$ is called
its {\it combinatorial type\/}.  We will sometimes blur
the distinction between a triangulation and its combinatorial
type.
There are $7$ different
combinatorial types of $8$-vertex triangulations of
tori. Figure~\ref{fig:triangulations} below shows the $8$-vertex triangulations
${\cal T\/}(46666677)$ and ${\cal T\/}(66666666)$.   We are
naming these particular triangulations by their degree sequence.
In each of these two cases, the degree sequence determines the
combinatorial type.

\begin{center}
  \resizebox{!}{2.6in}{\includegraphics{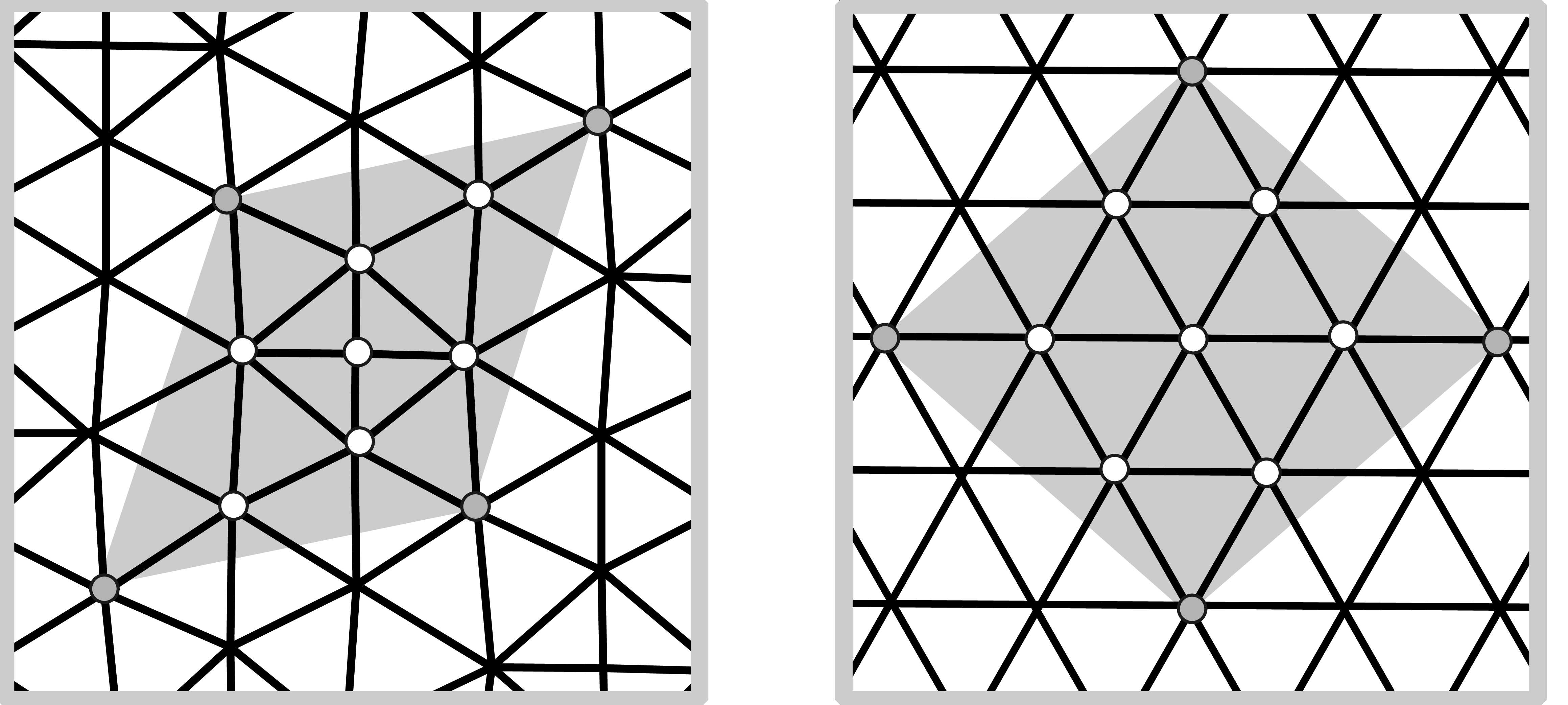}}
\figcaption{fig:triangulations}{The triangulations
  ${\cal T\/}(46666677)$ and ${\cal T\/}(66666666)$.}
\end{center}

An {\it origami torus\/} is an embedded topological
torus in $\R^3$ that is made by fitting together finitely
many triangles so that the cone angle around each vertex is $2\pi$.
Equivalently,
an origami torus is the image of a piecewise
affine isometric embedding
$\phi: T \to \R^3$
of a flat torus $T$.  This is to say that
$T$ has a triangulation  ${\cal T\/}$ such that
the restriction of $\phi$ to each triangle
of ${\cal T\/}$ is an isometric affine map.  (The specific
triangulation is part of the data.)
We call a lattice $\Lambda$ {\it realized\/} by ${\cal T\/}$ if there
exists an origami torus associated to the lattice $\Lambda$ and the
combinatorial type ${\cal T\/}$.  We call a combinatorial type
${\cal T\/}$ {\it universal\/} if every lattice $\Lambda$ is
realized by ${\cal T\/}$.  Here is our main result.

\begin{theorem}[Universality]
  ${\cal T\/}(46666677)$ is universal.
\end{theorem}

The Universality Theorem says that there is a single $8$-vertex combinatorial type out of which you can build an origami torus of any shape you like.  Our construction yields $2$ continuous families of $8$-vertex origami tori based on ${\cal T\/}(46666677)$ which together realize all shapes. We do not know if a single continuous family suffices, but see Remark~\ref{rem:onefamily}. Here $8$ is as good as it gets: In  [{\bf S1\/}] it is proved that $7$-vertex origami tori do not exist. See also [{\bf S2\/}].   Figure 1.2 below shows a
representative example of one of our tori.

\subsection{History and Context}
\label{sec:history}

The Universality Theorem is the culmination of a series
of results about origami tori, some quite recent.
\newline \newline
\noindent $\bullet$ In 1960, Y. Burago and V. Zalgaller 
[{\bf BZ1\/}] give the first construction
of origami tori.
\newline \newline
$\bullet$ In 1995, Y. Burago and V. Zalgaller  [{\bf BZ2\/}] prove that
every lattice is realized by some combinatorial type.
\newline \newline
\noindent $\bullet$ In 2020, T. Tsuboi [{\bf T\/}] gives
explicit constructions which realize every lattice as
an origami torus called  a {\it diplotorus\/} (For the rectangular lattices
  one needs a modified construction.)
P. Arnoux, S. Leli\`evre, and
A. M\'alaga [{\bf ALM\/}] give an independent
proof of the same result.
The construction requires infinitely many combinatorial types.
Diplotori were
discovered by U. Brehm [{\bf Br\/}]
in 1978.  See also [{\bf Se\/}, \S 6].
\newline \newline
\noindent $\bullet$ In 2023, F. Lazarus and F. Tallerie [{\bf LT\/}] give a universal triangulation having $1217$ vertices. Their paper also has an excellent discussion of the various attempts made to realize all lattices as origami tori.
\newline \newline
\noindent $\bullet$ In 2025, R.E.S. [{\bf S1\/}] constructs an
$8$-vertex origami torus and proves that
a $7$-vertex origami torus does not exist.
See also [{\bf S2\/}].  The $8$-vertex origami
tori in [{\bf S1\/}] and [{\bf S2\/}], called
{\it pup tents\/},
 have combinatorial type ${\cal T\/}(66666666)$ and
 are symmetric with respect to order $2$ rotation about the $Z$-axis.
 \newline \newline
 \noindent $\bullet$ In 2025,  P.G.D. and R.E.S.  [{\bf DS\/}] prove that a
 lattice is realized by a pup tent 
provided that it is neither rectangular nor sharp rhombic.
A lattice is {\it rectangular\/} if it is generated by a pair of perpendicular
vectors and {\it sharp rhombic\/} if it is generated by a pair of
equal-length vectors which make an angle of at most $\pi/3$ with
each other.  When the angle is $\pi/3$ you get the famous
{\it hexagonal lattice\/}.
Thus, almost every lattice is realized by a pup tent.
In [{\bf DS\/}] the authors speculate that perhaps the remaining
lattices can be realized by $8$-vertex origami tori without
$2$-fold symmetry.
\newline \newline
\noindent $\bullet$ In 2026, F. Lander  [{\bf L\/}] proves that the square lattice
is realized by ${\cal T\/}(66666666)$ and the hex
lattice is realized by ${\cal T\/}(46666677)$. These origami tori have
no symmetries.
\newline

The result in [{\bf L\/}] grew out
of some experiments that F.L. and S.T. did in the summer of
$2026$.  These authors made a poster for R.E.S.'s birthday
conference and, in an exciting development, their software
could realize some rectangular and sharp rhombic lattices.
Upon further experimentation, F.L. and S.T. discovered that
the folding pattern used to
realize the hexagonal lattice from ${\cal T\/}(46666677)$ turned
out to be extremely robust. This fact led to the Universality
Theorem.

\subsection{Proof Strategy and Organization}

Writing down explicit origami tori seems quite difficult, partly because
there are too many choices to make. The idea that emerged in [{\bf DS\/}], and was
used again in [{\bf L\/}], is to work instead on the boundary of the space of origami tori, where faces
collapse together and so make some of these choices for us. The most
collapsed examples lie entirely in a plane; we find these first and then
push them back inside. Here we carry this out for ${\cal T\/}(46666677)$
and every lattice. In \S\ref{sec:plane} we give a geometric account of the strategy,
ending in a proposition that reduces the proof to finitely many checks.
We verify these checks in three steps.

\medskip
\noindent{\bf Step 1:\/} We give an
explicit rational parametrization of a 5-dimensional family of {\it flattened
origami tori\/}: flat tori which lie entirely in the plane! Within it we find
two 2-parameter families which together realize every modulus. We do this
in \S\ref{sec:family}.

\medskip
\noindent{\bf Step 2:\/} For each family we give a single vertical lift,
the same heights for every member, which is embedded throughout. We call this {\it inflating\/}, and we do
it in \S\ref{sec:lift}.

\medskip
\noindent{\bf Step 3:\/} We check that horizontal motions can correct the
cone angles and the modulus, so that Proposition~\ref{prop:main} promotes
each (inflated) flattened origami torus to an origami torus of the same modulus.  We do this in
in \S\ref{sec:universality}.

\medskip
The two families came out of extending the searches performed in
[{\bf L\/}]. On several occasions we need certain polynomials not to vanish
on certain domains. To keep the story clean, we leave the description of
these certification techniques to \S\ref{sec:checks}.

Figure 1.2 shows the three step process which produces our
examples.

\begin{center}
\includegraphics[width=\textwidth]{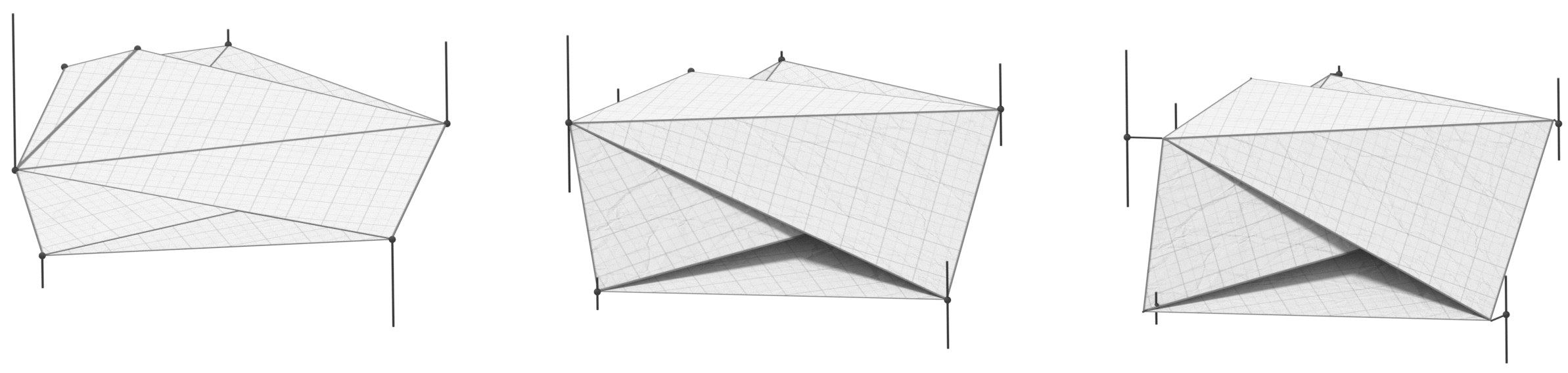}
\figcaption{fig:examples}{A flattened origami torus,
inflated then corrected}
\end{center}

  \subsection{Computer Calculations}

  We perform our calculations in Mathematica [{\bf W\/}].
  These calculations are all fairly short and easy to
  reproduce.  One file contains everything. Here are
  two places to get the file. \newline
  {\footnotesize
    {\bf http://www.math.brown.edu/$\sim$res/Java/UniversalityProof\/}
  \newline
  {\bf
    https://github.com/FabianLander/universal-8-vertex-origami-torus-triangulation\/}
  \/}

  \subsection{Artificial Intelligence Use}
  
  We used Artificial Intelligence in a moderate
  way.  The conceptual framework for the proof is
  ours, and evolved from the papers [{\bf S1\/}],
  [{\bf DS\/}], and [{\bf L\/}].
  At various stages, ChatGPT 5.6 implemented
  some of our search algorithms and optimized some of our code.
  Perhaps needless to say, the AI also helped format equations,
  provide references, and proofread the paper.
  Also, we mention that ChatGPT had an excellent idea for
  simplifying the angle-derivative formulas for
  the Jacobian calculation in \S\ref{sec:universality}.
  
  \subsection{Acknowledgements}

  R.E.S. would like to thank the I.H.E.S. for their support, in the form of the Gretchen and Barry Mazur Chair, during the writing of this paper.

\section{Producing Origami Tori}
\label{sec:plane}

We now describe the method we use to produce origami tori, following
[{\bf DS\/}] and [{\bf L\/}]: we construct degenerate example and then inflate them.  In this paper, our degenerate examples are ``completely flattened'': they lie in the plane.

\subsection{Looking on the Boundary}
\label{sec:boundary}

Fix the triangulation ${\cal T\/}={\cal T\/}(46666677)$, though the idea is
general. A polyhedral torus of type ${\cal T\/}$ is determined by its 8
vertices in $\R^3$, so these tori form the open subset of $\R^{24}$ where no
triangle is degenerate. The origami tori are its flat, embedded members,
${\cal O\/}={\cal E\/}\cap{\cal F\/}$, where ${\cal E\/}$ is the open set of embedded
tori and ${\cal F\/}$ is the set of intrinsically flat tori, whose triangle
angles sum to $2\pi$ at each vertex.

Origami tori are quite hard to produce: although 8-vertex origami tori are
known to exist [{\bf S1\/}], [{\bf DS\/}], [{\bf L\/}], none is known in closed
form. The two sets of constraints are hard in different ways. Embeddedness
is an open condition, so embedded tori can at least be found by sampling,
though they are rare: for our triangulation they make up approximately
$0.009\%$ of all tori.\footnote{From $10^9$ samples with respect to the
uniform measure on similarity classes of configurations: the round measure
on the unit sphere of configurations with centroid at the origin.}
Flatness is a system of equations: every cone angle must be exactly $2\pi$.
Sampling will never find a flat torus; they must instead be
constructed by solving a nonlinear system exactly, which means making many
choices.
One way forward is to impose extra conditions that make some of these
choices for us. But the resulting tori must still land in ${\cal E\/}$, which
is small. The idea is to draw the extra conditions from ${\cal E\/}$
itself.

Turning some of the inequalities of embeddedness into equalities puts us on
the boundary of ${\cal E\/}$. This boundary is stratified by how embeddedness
fails: a vertex can touch a face, an edge can lie in a face, faces can
coincide, and each further contact imposes more conditions. Carried to the
end, this suggests a natural family: tori whose triangles have all collapsed
into a single plane, which we take to be $\R^2\times\{0\}$. We call these
{\it planar tori\/}, write ${\cal P\/}$ for the space of them, and call the
flat ones {\it flattened origami tori\/}.

Planarity gives us enough structure to hold on to. In space, a triangle can rotate about an edge through a whole circle of dihedral angles; in the plane, it is either folded over or not. Once we choose which triangles are folded over, the only remaining condition is that the folding closes up around each vertex and each loop of the torus. These are easy equations in the Euclidean group of the plane, and in \S\ref{sec:family} we give explicit solutions for ${\cal T\/}$ across all of moduli space.

\subsection{Inflating from the Boundary}

We took the equations for planar tori from the boundary $\partial {\cal E\/}$, but
some solutions may not lie in $\partial \cal E$.
A simple sufficient condition
for a planar torus $\Psi$ to lie in $\partial \cal E$ is this:
If some vertical lift $\Psi+z$, giving the vertices
heights $z\in\R^8$, is embedded, then so are the scaled lifts $\Psi+tz$,
which converge to $\Psi$ as $t\to0$.

If $\Psi$ is also flat, it sits in ${\cal F\/}$ right at the edge of
${\cal E\/}$, and origami tori should lie just inside. To reach them, we push
$\Psi$ a little way into ${\cal E\/}$, to $\Psi+tz$. This leaves ${\cal F\/}$, so
we correct by moving the vertices horizontally, parallel to the plane, by
some $h\in\R^{16}$ that brings us back to ${\cal F\/}$.  The following
result assures us that the
correction does not push us back out of ${\cal E\/}$.

\begin{lemma}
\label{lem:vertical}
Let $\Psi$ be a planar torus such
that $\Psi+z$ is embedded. Then there is $\eta>0$ such that $\Psi+tz+h$ is
embedded for every $t>0$ and every horizontal motion $h$ with
$|h|<\eta$.
\end{lemma}

\startproof
The vertical scaling $(X,Y,Z)\mapsto(X,Y,Z/t)$ is an affine
bijection of $\R^3$, and it sends $\Psi+tz+h$ to $\Psi+z+h$. Since
embeddedness is open, $\Psi+z+h$ is embedded once $|h|$ is small enough.
\endproof

This is a lot of room: as $t\to0$ the lifted torus collapses onto the plane,
yet it tolerates horizontal motions of a fixed size. So any horizontal
correction that tends to $0$ with $t$ eventually fits inside it. For
flatness, the implicit function theorem supplies such a correction whenever
horizontal motions can adjust the cone angles at $\Psi$. But nothing in this
argument is special to cone angles. The same reasoning inflates $\Psi$ while
keeping fixed any quantity that depends smoothly on the vertices, provided
horizontal motions can adjust it at $\Psi$.

\begin{center}
\includegraphics[width=0.66\textwidth]{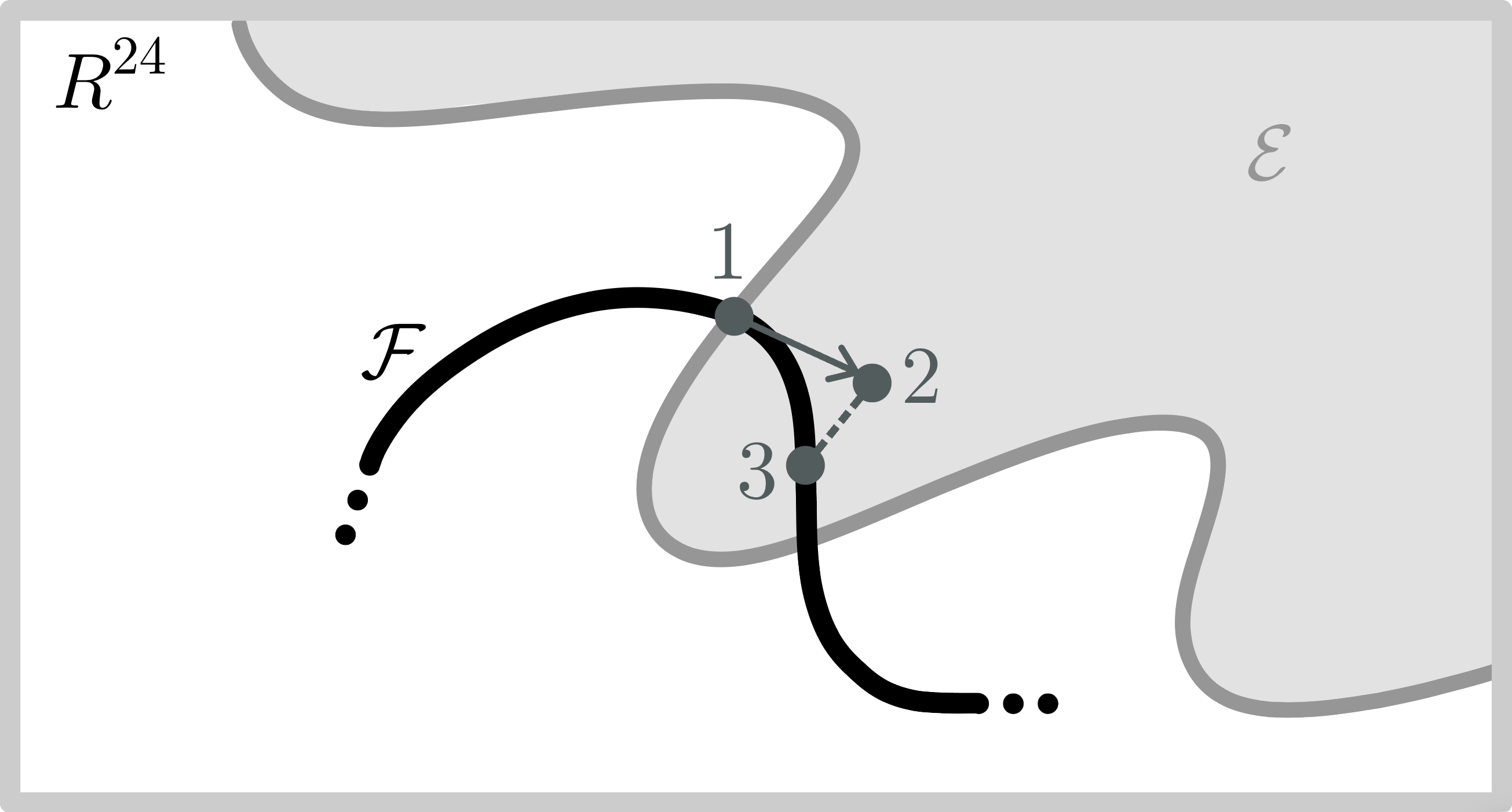}
\figcaption{fig:inflate}{Inflating from the boundary.}
\end{center}

\begin{proposition}
\label{prop:main}
Let $\Psi$ be a planar torus and $\Phi$ a smooth
map from a neighborhood of $\Psi$ to $\R^m$. Suppose that
\begin{enumerate}
\item\label{hyp:lift} $\Psi$ has an embedded vertical lift $\Psi+z$;
\item\label{hyp:submersion} $\Phi$ restricted to ${\cal P\/}$ is a submersion at $\Psi$.
\end{enumerate}
Then there are embedded tori arbitrarily close to $\Psi$ with the same value
of $\Phi$.
\end{proposition}

\startproof
By hypothesis~\ref{hyp:submersion} there are horizontal vectors
$e_1,\ldots,e_m$ with $d\Phi_\Psi(e_i)$ the standard basis of $\R^m$. Let
\[
G(t,\epsilon)=\Phi(\Psi+tz+\epsilon_1e_1+\cdots+\epsilon_me_m).
\]
Then $G(0,0)=\Phi(\Psi)$ and $\partial G/\partial\epsilon$ is the identity at
$(0,0)$, so the implicit function theorem gives a smooth $\epsilon(t)$ with
$\epsilon(0)=0$ and $G(t,\epsilon(t))=\Phi(\Psi)$ for small $t$. The
horizontal motion $h(t)=\sum\epsilon_i(t)e_i$ tends to $0$ with $t$, so by
Lemma~\ref{lem:vertical} the torus $\Psi+tz+h(t)$ is embedded for small $t>0$, and it has
the same value of $\Phi$ as $\Psi$.
\endproof

Proposition~\ref{prop:main} provides easily checkable conditions that let us promote a
planar torus with the properties we want to an inflated torus with the same
properties.

In this paper we want origami tori whose modulus we know exactly. So we must
control both the cone angles and the modulus, and in \S\ref{sec:universality} we take $\Phi$ to
record both. This makes Proposition~\ref{prop:main} the engine of the paper, and the
rest of the paper checks its hypotheses.
\section{The Flattened Tori}
\label{sec:family}

\subsection{The Families}
\label{sec:need}

Universality requires an origami torus for every lattice. Scaling in $\R^3$
takes origami tori to origami tori, and a lattice $\Lambda$ and its mirror
image $\overline\Lambda$ give isometric tori, so it is enough to realize
one lattice in each class up to scaling and reflection. Every such class
contains $\Lambda_\zeta=\Z+\Z\zeta$ for some $\zeta$ in
\[
{\cal M\/}=\{x+iy:\ 0\le x\le\tfrac12,\ x^2+y^2\ge1\},
\]
the closed $(2,3,\infty)$ triangle in ${\bf H}^2$, which is half of the
standard fundamental domain for the modular group. So for each
$\zeta\in{\cal M\/}$ we need a flattened origami torus with lattice
$\Lambda_\zeta$ satisfying the hypotheses of Proposition~\ref{prop:main}.

Since ${\cal M\/}$ is a continuum, we work with families of flattened tori
$\Psi(\zeta)$ with lattice $\Lambda_\zeta$, depending smoothly on $\zeta$;
along a family, each hypothesis of Proposition~\ref{prop:main} becomes a
statement about functions of $\zeta$, checked over a whole region at once.
A few families whose regions cover ${\cal M\/}$ then suffice.

Our families come from fixing how the torus folds. When a torus is folded
into the plane, each triangle lands either face up or face down, so every
planar torus has a {\it folding pattern\/}, recording which triangles are
reversed. For a fixed pattern, the flattened origami tori are found by
solving the equations that make the folding close up. For the
Universality Theorem we need only a single folding pattern, shown in
Figure~\ref{fig:square}. We spoil the story and give the formulas right
away. They form a $5$-parameter family of flattened origami tori, with
parameters $\zeta=x+iy$ and $u,k,\lambda$. Fixing $u,k,\lambda$ as
functions of $\zeta$ in two ways then gives two families, which together
cover ${\cal M\/}$.

Each member has two planar pictures: the flat torus before folding, described by vertices $w_j$ and their lattice translates, and its folded image, described by vertices $p_j$. The folding maps corresponding triangles isometrically. We begin with the folded vertices.

\paragraph{The Flattened Origami Torus:}
Here are the points of the corresponding flattened origami torus $\Psi$.
With the auxiliary variable $b=1-x-2k$, they are
\[
\begin{aligned}
p_0&=(0,\ 0), & p_1&=\Bigl(u,\ \frac{u(b+2k/\lambda)}{y}\Bigr),\\
p_2&=\Bigl(u-\tfrac12,\ \frac{(1-2u)x}{2y}\Bigr), & p_3&=\Bigl(-k-b\lambda,\ \frac{\lambda(y^2-b^2)-2bk}{2y}\Bigr),\\
p_4&=\Bigl(x+k-\tfrac12,\ \frac{y^2+bx}{2y}\Bigr), & p_5&=\Bigl(u-\tfrac12,\ \frac{y^2+b(b+2u-1)}{2y}\Bigr),\\
p_6&=\Bigl(u-k,\ \frac{b(u-k)}{y}\Bigr), & p_7&=\Bigl(u,\ \frac{y^2+b^2+2bu}{2y}-\frac{2u-2x-2k+1}{2y(1-\lambda)}\Bigr).
\end{aligned}
\]

\paragraph{The Flat Torus:}
We will describe $16$ triangles.  We get a tiling of $\R^2$ by taking the orbit of these triangles under the action of the lattice $\Lambda_\zeta$ defined above.  The quotient then gives us our $8$-vertex triangulation of the flat torus $T=\R^2/\Lambda_\zeta$.  Here are the triangles
\begin{equation}
  \label{pattern}
\begin{array}{cl@{\qquad}cl}
012 & (w_0,w_1,w_2) &
167 & (w_1,w_6,w_7)\\
031 & (w_0,w_3,w_1) &
264 & (w_2+e,w_6,w_4+e-\tau)\\
024 & (w_0,w_2,w_4-\tau) &
257 & (w_2+e,w_5+e,w_7)\\
043 & (w_0,w_4-\tau,w_3) &
276 & (w_2+e,w_7,w_6)\\
152 & (w_1,w_5,w_2) &
347 & (w_3+\tau,w_4,w_7)\\
136 & (w_1,w_3,w_6) &
356 & (w_3,w_5+e-\tau,w_6)\\
145 & (w_1,w_4,w_5) &
375 & (w_3+\tau,w_7,w_5+e)\\
174 & (w_1,w_7,w_4) &
465 & (w_4+e-\tau,w_6,w_5+e-\tau)
\end{array}
\end{equation}
Here $e=(1,0)$ and $\tau=(x,y)$. We have $w_j=p_j$ for $j\ne3,6$, and
\[
w_3=\Bigl(k,\ -\frac{2bk+\lambda(b^2+y^2)}{2y}\Bigr),\qquad
w_6=\Bigl(u+k,\ \frac{b(u-k)}{y}\Bigr).
\]

\paragraph{Preliminary Isometric Mapping:}
Let $W_j$ denote the $j$th triangle above. We call the set $\{W_0,...,W_{15}\}$ the {\it fundamental triangles\/}. Let $P_j$ denote the corresponding triangle whose vertices are the $p$ points. Let $\phi: \bigcup W_j \to \bigcup P_j$ be the piecewise affine map from the one set of $16$ triangles to the other.

\begin{lemma}
  \label{iso0}
$\phi$ is an isometry on each of the $16$ triangles.
\end{lemma}

\startproof
For each of the $16$ triangles, one checks that the three squared edge
lengths of $W_j$ and of $P_j$ agree; these are rational functions of
$x,y,u,k,\lambda$, so this is a direct symbolic computation. We do this
exactly with the file {\bf UniversalityProof\/}.
\endproof

\begin{remark}
\label{rem:constraints}
For $\zeta\in{\cal M\/}$, there are some algebraic constraints on the fiber $(u,k,\lambda)$ which suffice to ensure  that the triangles above have positive area:
\begin{equation}
u,\ k,\ x+k-u,\ \lambda/2\in(0,1/2),\hskip 30 pt
\frac{u-x-k+1/2}{1-\lambda}+\frac{2ku}{\lambda}
<b+2k(x+k).
\end{equation}
We will never need to check or work with these constraints.
\end{remark}

\begin{remark}
\label{rem:onefamily}
By the constraint $x+k-u\in(0,\frac12)$, $k>u$ when $x=0$ and $k<u$ when
$x=\frac12$. So a continuous family over all of ${\cal M\/}$ must pass
through $k=u$, where $p_0=p_6$, and there it seems impossible to add
heights that give an embedded lift. This is why we cover ${\cal M\/}$
with two families.
\end{remark}

\paragraph{The Two Families:}
We now pick out two subfamilies, parametrized by the modulus $\zeta=x+iy$, whose domains together cover ${\cal M\/}$:
\begin{alignat*}{2}
&\text{Domain 1}=\{\zeta\in{\cal M\/}: x\le 1/6\},\quad
&&\text{Family 1: }(u,k,\lambda)=(3/128,\,11/32,\,13/32),\\
&\text{Domain 2}=\{\zeta\in{\cal M\/}: x\ge 1/6\},\quad
&&\text{Family 2: }(u,k,\lambda)=(x/2-1/2048,\,1/16,\,1/2).
\end{alignat*}
Both formulas make sense for $\zeta$ in an open neighborhood of their domains: the denominators never vanish there.

These choices of family give the vertices above simple coordinates.
Every first coordinate is linear in $x$, and every second coordinate is
a rational function $(f(x)+c\,y^2)/y$, where $f$ is at most quadratic
and $c$ is constant. This greatly simplifies the polynomials in our checks: the
triangle areas (Lemma~\ref{lem:fam}), the volumes of \S\ref{sec:lift}
(Lemma~\ref{vol}), and the Jacobian (Lemma~\ref{subd}).

In \S\ref{sec:produce} we check that the members of the two families are
flattened origami tori; this is all the proof needs. In
\S\ref{sec:finding} we show how one can derive these formulas from the
folding pattern.

\subsection{Checking the Families Are Flattened Origami Tori}
\label{sec:produce}

\begin{center}
\begin{minipage}[c]{0.4\textwidth}
  \centering
  \includegraphics[width=\textwidth]{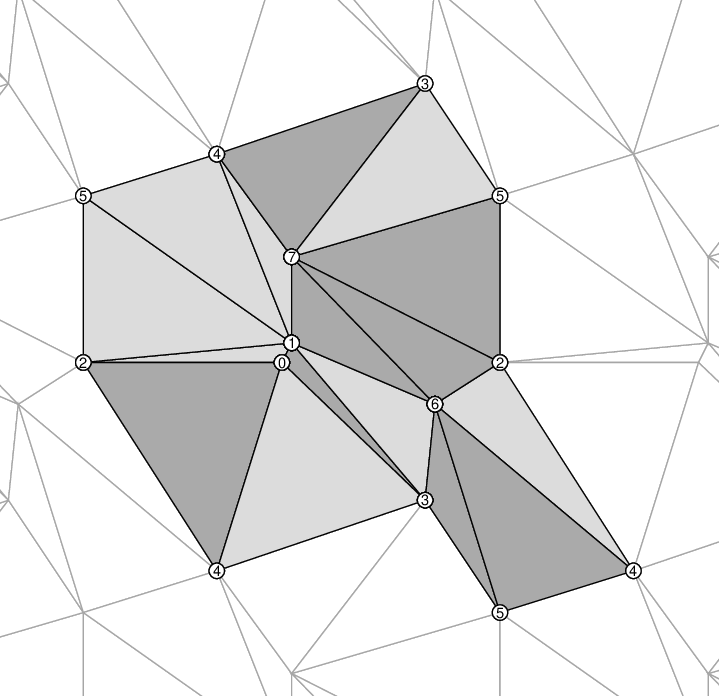}
\end{minipage}\hfill
\begin{minipage}[c]{0.55\textwidth}
  \centering
  \includegraphics[width=\textwidth]{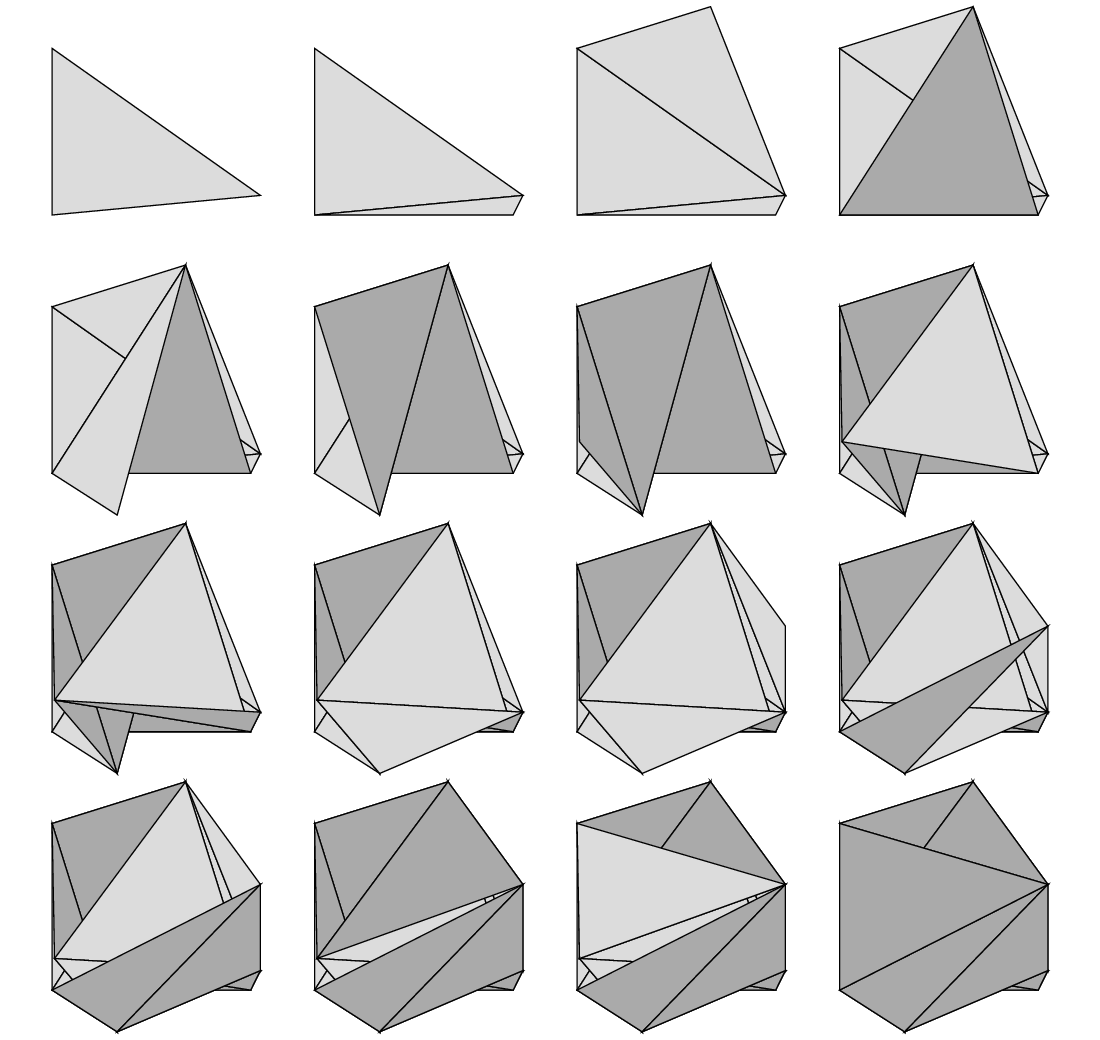}
\end{minipage}
\figcaption{fig:square}{The square torus and its flattened origami realization.}
\end{center}

\begin{center}
\begin{minipage}[c]{0.4\textwidth}
  \centering
  \includegraphics[width=\textwidth]{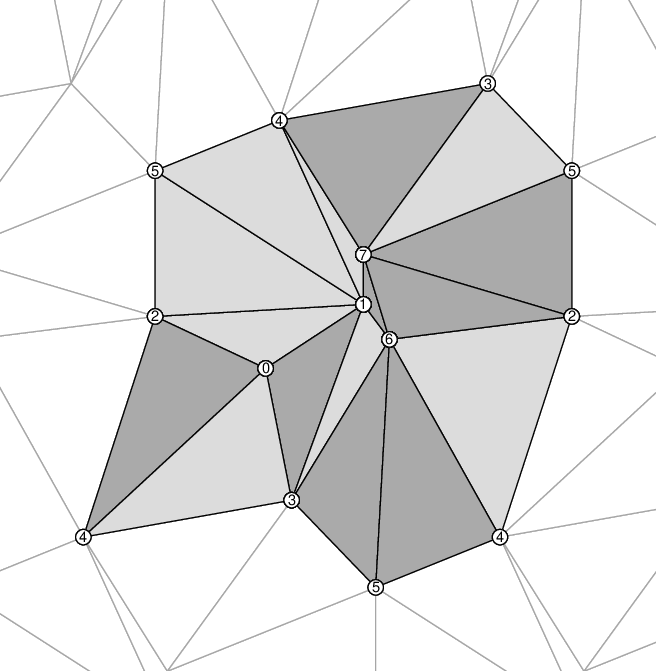}
\end{minipage}\hfill
\begin{minipage}[c]{0.55\textwidth}
  \centering
  \includegraphics[width=\textwidth]{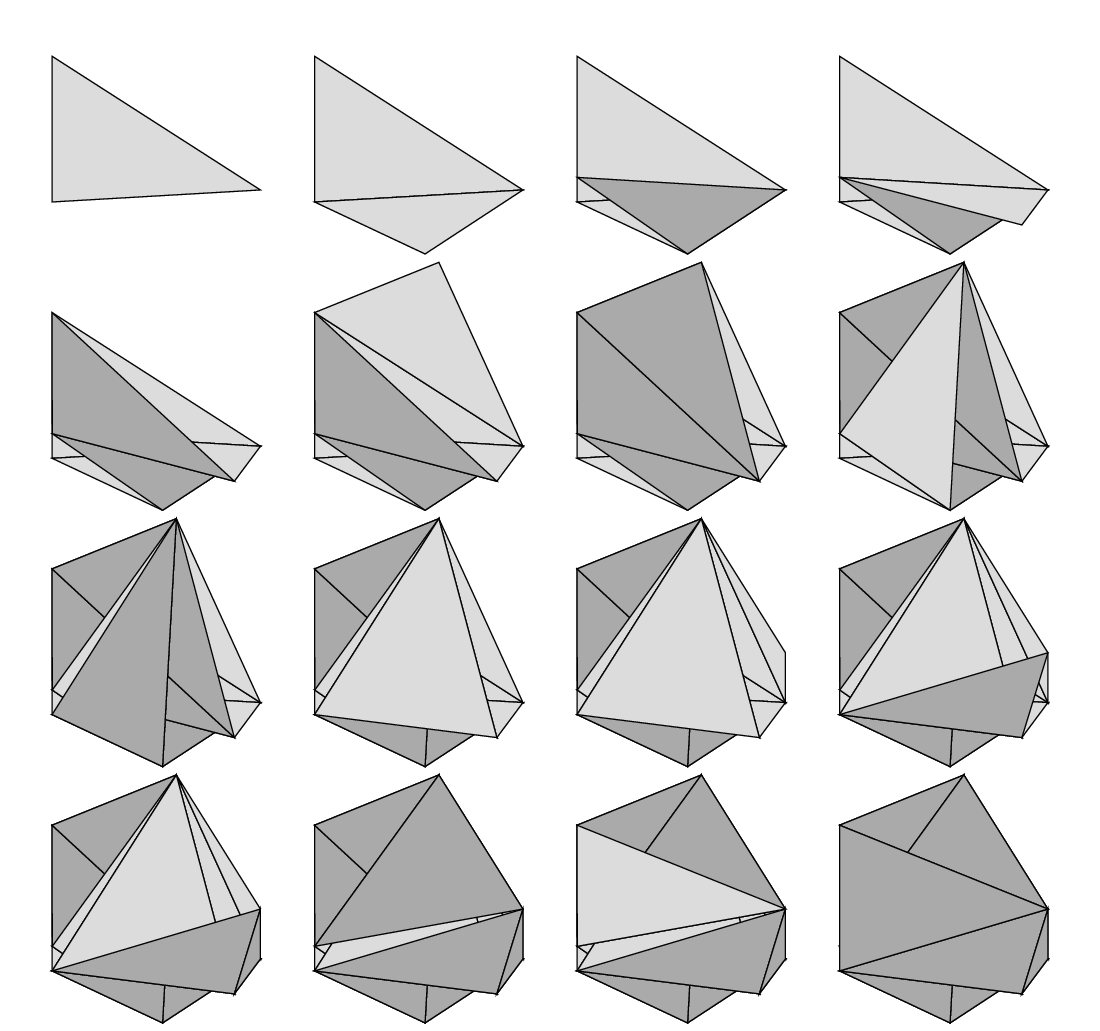}
\end{minipage}
\figcaption{fig:family2}{The torus of modulus $8/17+i$ and
  its flattened  origami realization.}
\end{center}

The left side of Figure~\ref{fig:square} shows the member of Family 1 at $\zeta=i$.  The right side shows the flattened origami torus built one triangle at a time, in an ordering that is compatible with our corresponding lift.  Figure~\ref{fig:family2} does the same thing for the somewhat randomly chosen member $\zeta=8/17+i$ of Family 2.  (It would be more mathematically satisfying to show the picture for the hex torus, but the picture does not look that good: two of the vertices are quite close together.)  The coloring convention in both figures is as follows. Fixing an orientation of $T$, light triangles keep their orientation and dark ones are reversed. Read down the left column of the table in \S\ref{sec:need} and then the right, with $+$ for kept and $-$ for reversed, our folding pattern is $+--+++++-+----+-$.

In both figures, the sixteen triangles form a polygon whose boundary
edges are paired by the lattice: $e$ pairs the edges $42$ and $25$ on the
left with those on the right, $\tau$ pairs $43$ and $35$ on the bottom
with those on the top, and $e-\tau$ pairs the two edges $45$. The table
shows that these pairings hold for every member, not only the two
pictured: for instance, triangles $145$ and $465$ share the edge $45$,
shifted by $e-\tau$.

\begin{lemma}
\label{lem:fam}
Let $d=1,2$.
\begin{enumerate}
\item\label{fam:area} For all parameters in a suitably chosen open neighborhood of Domain $d$, every triangle $W_j$ in Family $d$ has positive area.
\item\label{fam:torus} Throughout a connected open neighborhood $U_d$ of
  Domain $d$, $\phi$ gives a flattened origami torus $\phi:T\to\Psi$.
\end{enumerate}
\end{lemma}

\startproof
(\ref{fam:area}) Twice the signed area of a triangle is the $2\times2$
determinant of two of its edge vectors, so we must show that these
sixteen determinants do not vanish. {\bf UniversalityProof\/} certifies
this exactly via strip dominance (\S\ref{sec:checks}), on each family's
domain.

One may also check this by hand. As the vertex coordinates are simple
(\S\ref{sec:need}), so are these determinants: each has the form
$A(x)/y+B(x)\,y$, so $y$ times it is $A(x)+B(x)\,y^2$, linear in $y^2$.
Here Lemma~\ref{lem:nonvanishing} reduces each triangle to a pair of sign
checks, for two at most quadratic polynomials on the family's domain, and
all of them hold. This open condition persists in an open neighborhood
of Domain $d$.

(\ref{fam:torus}) Let $\Pi=\bigcup W_j$ and $\Psi=\bigcup P_j$. We have a map $\phi:\Pi\to\Psi$ which is an isometry on each triangle, by Lemma~\ref{iso0}. For each of the figures we plotted above, $\Pi$ is an embedded solid polygon and the three elements $e,\tau,e-\tau$ of $\Lambda_\zeta$ pair the sides of $\Pi$. Let $T$ be the quotient of $\Pi$ by this side-pairing. $T$ is a topological torus. Because the side-pairings are translations, we have the familiar result that $T$ is a flat torus. This picture persists for an open set of parameters in Domain $d$, for $d=1,2$.

We let $U_d$ be some connected open neighborhood of Domain $d$ for which (\ref{fam:area}) holds. Then (\ref{fam:area}) tells us that the triangles retain their orientations and never degenerate in $U_d$. Also, all the angles vary analytically. Hence, if we form the abstract topological space by gluing together the triangles in the same pattern as we do when we know that $\Pi$ is embedded, we get a topological torus whose cone angles vary analytically. Since these cone angles equal $2\pi$ on an open set, they are always $2\pi$. Hence, our torus $T$ is always a flat torus. (The expert reader may enjoy proving the unnecessary fact that $\Pi$ is always an embedded polygon and a fundamental domain for $\Lambda_\zeta$.) Furthermore, $\phi$ respects all the gluings needed to make $T$ and therefore gives us a flattened origami torus $\phi:T\to\Psi$ throughout $U_d$ for each choice of $d$.
\endproof

\subsection{Producing Families by Folding}
\label{sec:finding}

In this subsection, which is not needed for the proof, we derive the
formulas in \S\ref{sec:need} from our folding pattern. There are very many
folding patterns, $2^{16}$ for each triangulation, and ours has several
properties that let the proof go through; we explain how we found it in
\S\ref{sec:found}.

Let $\Psi=\phi(T)$ be any flattened origami torus with our folding
pattern, where $\phi$ is isometric on each triangle of its intrinsic
flat torus $T$. An edge is a {\it crease\/} if its two triangles have
opposite signs. Two triangles of the same sign sharing an edge are moved
by the same isometry. So if we cut the universal cover $\R^2$ of $T$
along the creases into {\it islands\/}, then $\phi$, lifted to $\R^2$,
moves each island by a single isometry, and the isometries of two
islands meeting along a crease differ by the reflection in it. For our
folding pattern the creases cut $\R^2$ into ten islands $A,\ldots,J$, up
to translation by the lattice.   

Figure~\ref{fig:islands} below shows these islands, drawn in the universal cover and shaded according to their orientation as in Figures 2.1 and 2.2.  More precisely, the folding keeps the orientation of the light triangles and reverses the orientation of the dark ones.  Our subtle shading in the figure is designed to distinguish the $16$ fundamental triangles from their images under the lattice action. We also note one additional bookkeeping feature of Figure 3.3:  A shaded vertex labeled $v$ is the translate of $w_v$ by the
indicated vector.

\begin{center}
\includegraphics[width=0.8\textwidth]{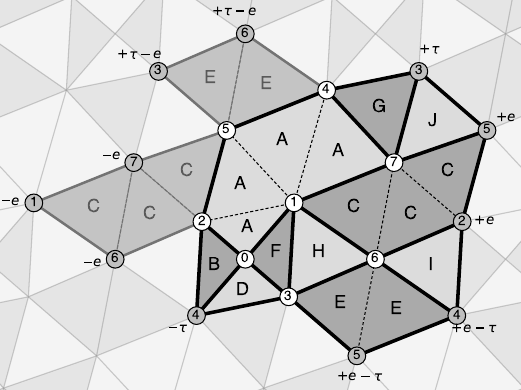}
\figcaption{fig:islands}{The islands and the sixteen triangles, drawn in
the universal cover.}
\end{center}

Our aim is to find all flattened origami tori with this folding pattern, up to similarity. Knowing how $\phi$ moves one island determines how it moves all the others, by crossing creases, so rather than start with $T$ and look for $\phi$, we build both together, one island at a time. Normalize the lattice of $T$ to be generated by $e=(1,0)$ and $\tau=(x,y)$ with $y>0$, and put $w_0=0$. Then $T$ is determined by lifts $w_j\in\R^2$ of its vertices, which are our unknowns. The table in \S\ref{sec:need} lists the sixteen triangles with the lifts of their vertices.

Composing with an isometry of the plane, we may take $\phi$ to be the
identity on the island $A$, and we write $p_j=\phi(w_j)$ for the vertices of
$\Psi$. Crossing creases then determines $\phi$ everywhere, and the only
requirement is that it be well defined on $T$: all copies of a vertex in
$\R^2$ must go to the same point. These {\it closing conditions\/} are the
only equations to solve, and they are easy. Each says that a composition of
reflections in creases sends one copy $P$ of a vertex to another copy $Q$.
When it is a single reflection, the crease is the perpendicular bisector of
$P$ and $Q$, so each point $c$ on it satisfies $|c-P|^2=|c-Q|^2$, which is
linear in $c$. For our folding pattern the conditions can be ordered so that
each is linear in the unknowns not yet determined. So we can solve them one
at a time, and the solutions form a family, rational in $x$, $y$ and the
choices left free.

We illustrate, through two examples, how the closing conditions are solved, starting from $A$. The first example shows both how a condition determines unknowns and how a free parameter arises. Since $\phi$ is the identity on $A$, $p_j=w_j$ for its vertices $0,1,2,4,5,7$. The island $C$ meets $A$ along the crease $17$, so $\phi$ moves $C$ by the reflection in $17$ (Figure~\ref{fig:folds}). $C$ contains the copies $w_2+e$ and $w_5+e$ of the vertices $2$ and $5$, while $A$ contains $w_2$ and $w_5$, so this reflection sends $w_2+e$ to $w_2$ and $w_5+e$ to $w_5$. Thus $17$ is the perpendicular bisector of $w_2$ and $w_2+e$, a vertical line $X=u$, and nothing determines $u$: it is our first free choice. It places $w_1$ and $w_7$ on $X=u$, and $w_2$ and $w_5$ on $X=u-\tfrac12$.

\begin{center}
\includegraphics[width=\textwidth]{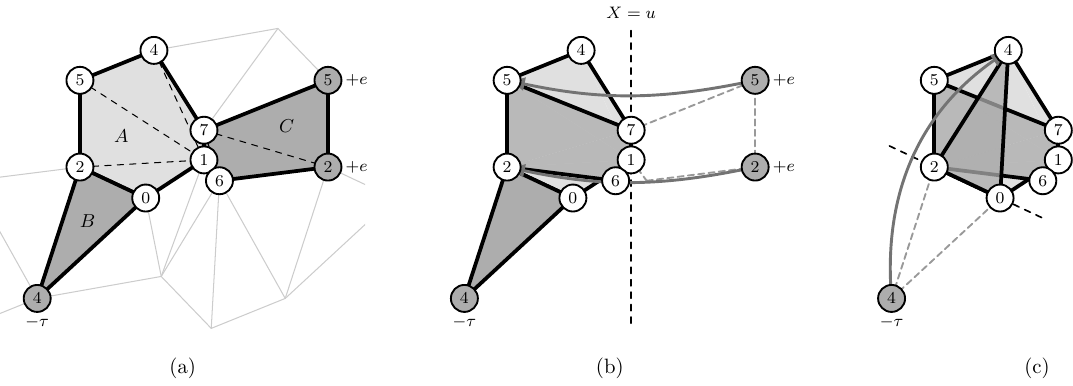}
\figcaption{fig:folds}{(a) $A$, $B$, $C$, unfolded.
(b) Fold $C$ across $17$. (c) Fold $B$ across $02$.}
\end{center}

Crossing into $B$ along the crease $02$ works the same way (Figure~\ref{fig:folds}c):
the reflection in $02$ must send $w_4-\tau$ to $w_4$, so $02$, which passes
through $w_0=0$, is perpendicular to $\tau$, and this pins down $w_2$ on
$X=u-\tfrac12$. The rest goes the same way once $w_3$ is placed freely,
which gives our other two choices: its first coordinate $k$ and its second
coordinate, which we parametrize by $\lambda$. Comparing $E$
with $D$ then fixes $w_4$ and $w_5$, comparing $E$ with $C$ fixes $w_6$, the
islands $F$ and $G$ fix $w_1$ and $w_7$, and $H$, $I$ and $J$ add no new
equations. Three free choices is what one expects: up to similarity, eight
points in the plane have 12 degrees of freedom, and flatness imposes 7
independent conditions, leaving the modulus and three more.

Working out all the closing conditions in this way led to the formulas of \S\ref{sec:need} for the vertices $p_j=\phi(w_j)$ of $\Psi$. The flat torus $T$ is given by the lifts $w_j$. Since $\phi$ is the identity on $A$, $w_j=p_j$ for the vertices $0,1,2,4,5,7$ of $A$. The other two vertices lie on islands moved by a single reflection: $w_3$ on $F$, reflected in the crease $01$, and $w_6$ on $C$, reflected in the crease $17$. This gives the formulas for $w_3$ and $w_6$.
\section{The Lifting Property}
\label{sec:lift}

Hypothesis~\ref{hyp:lift} of Proposition~\ref{prop:main} asks for an embedded vertical
lift. For each family we propose one set of heights
$z=(z_0,\ldots,z_7)$ for every member:
\[
\begin{array}{l@{\quad}*{8}{wc{2.8em}}}
 & z_0 & z_1 & z_2 & z_3 & z_4 & z_5 & z_6 & z_7\\[2pt] \hline
 \rule{0pt}{14pt}\text{Family 1} & 0 & 0 & -11 & -1 & -3 & -15 & 0 & 977\\
 \text{Family 2} & 0 & -2017 & 0 & -2 & -1 & -3 & 0 & 1430
\end{array}
\]
If $\Psi$ has vertices $p_j$, the vertical lift $\Psi+z$ has vertices
$\widehat p_j=(p_j,z_j)$. To see that it is embedded, we must check that no
two of its triangles meet except along a common edge or vertex. We do this
in two steps: we check that the lift is embedded at a single modulus, and
then that nothing changes across the family.

Throughout we use tetrahedra. Any four vertices span a tetrahedron, and
the sign of its volume records their relative position: on which side of
the plane through three of them the fourth one lies. We first show that
none of these $70$ volumes vanishes on the domain, so their signs are
constant across the family. We then explain how the signs certify
embeddedness.

\begin{lemma}
\label{vol}
For both families and every modulus in the domain, no four vertices of
the lift $\Psi+z$ are coplanar, that is, they are in general position.
\end{lemma}

\startproof
For each quadruple $\alpha,\beta,\gamma,\delta$ of vertices we must show that
$$V=\det(\widehat p_\alpha-\widehat p_\delta,\;\widehat p_\beta-\widehat p_\delta,\;\widehat p_\gamma-\widehat p_\delta)$$
does not vanish on Domain $d$. This is certified in {\bf UniversalityProof\/}
using the tweaked method of strip dominance from \S\ref{sec:checks}.

One may also check this by hand. The columns of $V$ are differences of
$X$-, $Y$- and $Z$-coordinates, so by \S\ref{sec:need}, multiplying the
middle column by $y$ and expanding along it gives
\[
yV=A(x)+B(x)\,y^2.
\]
For each of the $70$ determinants, Lemma~\ref{lem:nonvanishing} then
reduces nonvanishing to the signs of $A+B(1-x^2)$ and $B$, two
polynomials of degree at most two, on the family's domain. All these
checks hold.
\endproof

\begin{lemma}
\label{embed00}
Suppose that the vertices of a polyhedral torus are in general position,
and that every edge is disjoint from every triangle with which it shares
no vertices. Then the torus is embedded.
\end{lemma}

\startproof
We show that any two triangles $\Delta_1$ and $\Delta_2$ intersect
precisely in their common vertices and edge. If $\Delta_1$ and $\Delta_2$
share an edge, their planes are distinct by general position. Hence they
intersect exactly along their common edge. Now suppose $\Delta_1$ and
$\Delta_2$ do not share an edge. If they have an unwanted intersection,
some edge of one triangle must intersect the other triangle outside their
common vertices. An edge having a vertex in common with the other triangle
cannot do this: by general position, it intersects the plane of that
triangle only at the common vertex. Thus the unwanted intersection must
involve an edge and a triangle with no vertices in common. This
contradicts our hypothesis.
\endproof

Next we test whether a segment $(a_0a_1)$ is disjoint from a triangle
$(a_2a_3a_4)$. Deleting $a_k$ from the list leaves the vertices of a
tetrahedron. Let $V_k$ be its signed volume, with the points taken in
increasing order. We assume the five points are in general position.

\begin{lemma}
\label{win}
$(01)$ and $(234)$ are disjoint provided $\max(V_0V_1,V_2V_3,V_3V_4)>0$.
\end{lemma}

\startproof
Let $[abcd]$ denote the signed volume of the tetrahedron with vertices
listed in the order $a,b,c,d$. If $V_1=[0234]$ and $V_0=[1234]$ have the
same sign, then $(01)$ lies strictly on one side of the plane containing
$(234)$. This gives disjointness. If $(01)$ intersects $(234)$, then
$[0123]$, $[0134]$ and $[0142]$ have the same sign. Since $[0123]=V_4$,
$[0134]=V_2$ and $[0142]=-V_3$, this fails if $V_2V_3>0$ or $V_3V_4>0$.
\endproof

\begin{corollary}
\label{cor:embedded}
For both families and every modulus in the domain, the lift $\Psi+z$ is
embedded.
\end{corollary}

\startproof
There are $138$ pairs $(e,\Delta)$ where $e$ is an edge and $\Delta$ a
triangle of ${\cal T\/}$ with no common vertex. By Lemmas~\ref{vol}
and~\ref{embed00}, $\Psi+z$ is embedded if all of them satisfy the criterion
of Lemma~\ref{win}. Using the file {\bf UniversalityProof\/} we check this
at $\zeta=i$ in Family 1 and at $\zeta=8/17+i$ in Family 2
(Figures~\ref{fig:square} and~\ref{fig:family2}). Each $V_k$ in the
criterion is the volume of four lifted vertices, so by Lemma~\ref{vol} it
never vanishes on the domain. Since the domain is connected, its sign is
constant there, and the criterion holds throughout.
\endproof

\begin{remark}
Constant heights that work for a whole family are not automatic: we
searched for them alongside the choice of family. At a single modulus the
checks above take a computer almost no time, so candidate heights are easy
to test at many moduli. Even so, they are delicate: changing a single
height by $\pm1$ usually breaks the argument. See \S\ref{sec:found}.
\end{remark}
\section{Universality}
\label{sec:universality}

\subsection{The Map $\Phi$}

By \S\ref{sec:lift}, each of our flattened origami tori has an embedded vertical lift,
so it can be inflated. Proposition~\ref{prop:main} assures us we can keep
the values of {\it any\/} smooth map $\Phi$ fixed during inflation,
provided it is a submersion on the planar tori at our starting point. So we
build a $\Phi$ to suit our purpose, recording both the cone angles and the
modulus of polyhedral tori.

For a polyhedral torus, let $\theta_v$ be the cone angle at the vertex $v$,
the sum of the angles of the triangles meeting there, and let
$\delta_v=2\pi-\theta_v$ be its angle deficit. Each angle of a
nondegenerate triangle is a smooth function of its edge lengths, so the
$\delta_v$ are smooth functions on the space of polyhedral tori. By
Gauss--Bonnet, the eight deficits of a torus sum to zero. So a polyhedral
torus is flat exactly when $\delta_0,\ldots,\delta_6$ vanish, and we record
these seven deficits as $\delta=(\delta_0,\ldots,\delta_6)$.

There is one subtlety: the modulus is only defined for flat tori, but
Proposition~\ref{prop:main} needs $\Phi$ on a neighborhood of $\Psi(\zeta)$,
which contains tori that are not flat. The modulus of a flat torus is read
off from its development into the plane: the holonomy along two generating
loops is a pair of translations, and the modulus is their ratio. Other tori
are cone surfaces, and technically the developing map is only defined on the universal cover of the punctured surface we get by removing the cone points.   However, we can define a kind of pseudo modulus, sufficient for our purposes, by choosing (once and for all) flat annuli bounding each generating loop and then defining the holonomy relative to these annuli.

\begin{lemma}
\label{lem:pseudomodulus}
Developing along a fixed choice of annuli defines a smooth function
$\tilde\mu$ of the edge lengths, on an open set of polyhedral tori
containing ${\cal F\/}$, which equals the modulus on ${\cal F\/}$.
\end{lemma}

\startproof
  When we lift these annuli to the universal cover of ${\cal T\/}$ we get {\it strips\/} along which we can properly define the developing map.
In the universal cover, choose a vertex $v$, a triangle
$\Delta$ containing $v$, and their translates $v_1,\Delta_1$ by $e$ and
$v_2,\Delta_2$ by $\tau$. For each $j$, choose a strip of triangles from
$\Delta$ to $\Delta_j$, each sharing an edge with the next (Figure~\ref{fig:chains}
shows our choice).

Given a polyhedral torus, develop each strip into the plane. Place
$\Delta$ by a fixed rule, for instance with $v$ at the origin, a fixed edge
of $\Delta$ at $v$ along the positive $Y$-axis, and $\Delta$ on the side
given by the orientation of the torus. Then attach each subsequent triangle
along its shared edge, on the opposite side. Each placement is determined
by the edge lengths and depends smoothly on them, because the triangles are
nondegenerate. So, identifying $\R^2$ with $\C$,
\[
\tilde\mu=\frac{{\rm dev}_2(v_2)-{\rm dev}_2(v)}{{\rm dev}_1(v_1)-{\rm dev}_1(v)}
\]
is a smooth function of the edge lengths wherever the denominator is
nonzero.

When the torus is flat, the two strips are part of its development, so
the differences are its holonomy translations, up to a common rotation,
and $\tilde\mu$ is the modulus. In particular the denominator is nonzero
on ${\cal F\/}$, and so on an open set containing it, where $\tilde\mu$ is
the required function.
\endproof

\begin{center}
  \includegraphics[width=0.8\textwidth]{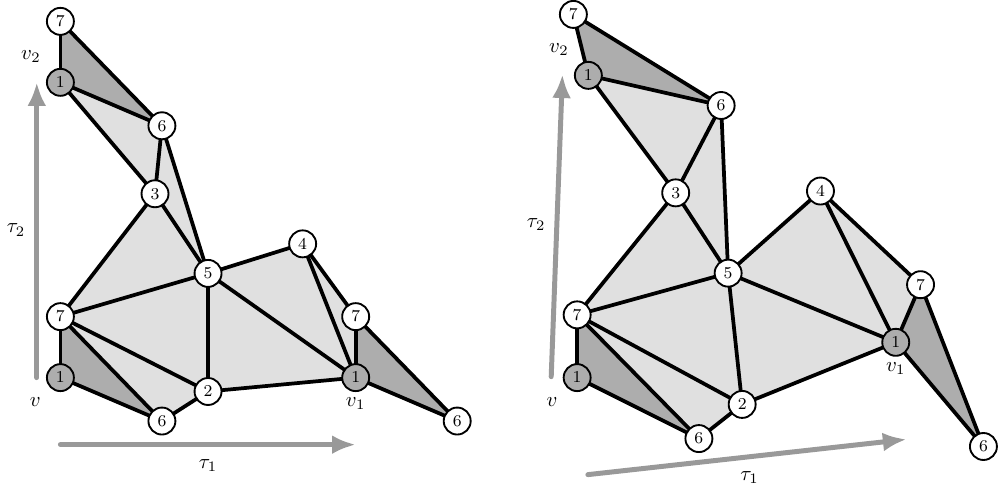}
\figcaption{fig:chains}{The modulus defined by strip development.}
\end{center}

Figure 5.1 shows our choice of strips, and also illustrates our construction of this pseudo modulus.
The strips for ${\cal T\/}(46666677)$, with $v$ at vertex $1$
  and $\Delta=167$ (dark), and the edge $17$ vertical: $$167,276,257,152,145,174,167$$
  from $\Delta$ to $\Delta_1$, and $$167,276,257,375,356,136,167$$ from $\Delta$
  to $\Delta_2$. They are developed for the
square torus (left) and for a nearby torus that is not flat (right).

Set $\Phi=(\tilde\mu,\delta)$, a smooth map from a neighborhood of each
$\Psi(\zeta)$ to $\C\times\R^7$. Near $\Psi(\zeta)$, the tori with
$\Phi=(\zeta,0)$ are exactly the flat tori of modulus $\zeta$: $\delta=0$
means flat, and on flat tori $\tilde\mu$ is the modulus.

To apply Proposition~\ref{prop:main} we must show that $\Phi$ restricted to
${\cal P\/}$ is a submersion at $\Psi(\zeta)$. Along our families the
modulus part is automatic, so this reduces to a statement about $\delta$
alone:

\begin{lemma}
\label{lem:reduction}
Let $\Psi(\xi)$ be a smooth family of flattened origami tori of modulus $\xi$, for $\xi$ in an open set $U\subset{\bf H}^2$ containing $\zeta$. If $\delta$ restricted to ${\cal P\/}$ is a submersion at $\Psi(\zeta)$, then so is $\Phi$ restricted to ${\cal P\/}$.
\end{lemma}

\startproof
Since each $\Psi(\xi)$ is planar, the velocities
$f_1=\partial\Psi/\partial x$ and $f_2=\partial\Psi/\partial y$ at $\zeta$
are horizontal. Since $\delta(\Psi(\xi))=0$ and $\tilde\mu(\Psi(\xi))=\xi$,
we have $$d\delta(f_j)=0, \hskip 20 pt
d\tilde\mu(f_1)=1, \hskip 20 pt \tilde\mu(f_2)=i.$$ By
hypothesis there are horizontal vectors $e_1,\ldots,e_7$ with
$d\delta(e_i)$ the standard basis of $\R^7$. On these nine vectors the
derivative of $\Phi$ has the block form
\[
\begin{pmatrix} I_2 & B \\ 0 & I_7 \end{pmatrix},
\]
which is invertible. So $d\Phi$ maps the horizontal vectors onto
$\C\times\R^7$.
\endproof

\subsection{Completing the Proof}

By Lemma~\ref{lem:reduction}, it remains to show that $\delta$ restricted to
${\cal P\/}$ is a submersion at each $\Psi(\zeta)$. In fact the first
coordinates of vertices $0,\ldots,6$ already suffice. Given
$\epsilon\in\R^7$, let $\Psi(\epsilon)$ move each vertex $p_i$ of $\Psi$
by $(\epsilon_i,0)$ for $i=0,\ldots,6$, leaving vertex $7$ fixed, and for
Family $d$ set
\[
J_d=\det\frac{\partial\delta}{\partial\epsilon}\bigg|_{\epsilon=0}.
\]
We show that $J_d$ does not vanish.

Writing $\zeta=x+iy$, $J_d$ is a rational function of $x$ and $y$, even
though the angles themselves are not. The angle $\alpha$ between vectors
$V=(v_1,v_2)$ and $W=(w_1,w_2)$ is
\[
\alpha=\arccos\bigg(\frac{V\cdot W}{\|V\|\,\|W\|}\bigg).
\]
This formula is far from rational, but its derivatives are, for example
\[
\frac{\partial\alpha}{\partial v_1}=\sigma\,\frac{v_2}{v_1^2+v_2^2},
\qquad
\sigma=\operatorname{sign}\det(V,W).
\]
Taking $V$ and $W$ in the order given by the orientation of the torus,
$\sigma$ is the sign of the triangle in our folding pattern. The vertex
coordinates of our flattened tori are rational functions of $x$ and $y$,
so each such derivative is too. The entries of
$\partial\delta/\partial\epsilon$ are sums of these derivatives, and $J_d$
is a polynomial in the entries, so $J_d$ is rational in $x$ and $y$. Its
denominators are squared edge lengths $v_1^2+v_2^2$, which are positive
because our triangles are nondegenerate (Lemma~\ref{lem:fam}), so $J_d$ is
defined on all of Domain $d$.

\begin{lemma}
\label{subd}
The function $J_d$ does not vanish on Domain $d$.
\end{lemma}

\startproof
This is certified in {\bf UniversalityProof\/} using the tweaked method of
strip dominance from \S\ref{sec:checks}. Alternatively, with a computer algebra system one
may verify that the numerator of $J_d$ in lowest terms factors as
$$C_d\,y^{13}\,b_d\,\bigl(A_d(x)+B_d(x)\,y^2\bigr),$$
where $C_d$ is a nonzero constant, $b_1=\tfrac5{16}-x$, $b_2=\tfrac78-x$,
and $A_d$ and $B_d$ are polynomials of degree $5$ and $3$. For example,
for Family 1,
{\small
\begin{align*}
A_1(x)&=13791139987456x^5-29311757123584x^4-347158265856x^3\\&\qquad+16581166028288x^2-5729028705984x+420099304545,\\
B_1(x)&=-13791139987456x^3+20692294631424x^2+16351578693632x+715752410880.
\end{align*}}
On Domain $d$, both $y$ and $b_d$ are positive. By Lemma~\ref{lem:nonvanishing}, the remaining factor is nonzero because the quintic $A_d+B_d(1-x^2)$ and the cubic $B_d$ are both positive for Family 1 and both negative for Family 2.
\endproof

This completes the proof of the Universality Theorem. Indeed, let
$\zeta\in{\cal M\/}$. Since our two domains cover ${\cal M\/}$, $\zeta$
lies in one of them, say Domain $d$, and $\Psi(\zeta)$, the flattened
origami torus of Family $d$ with modulus $\zeta$, satisfies both
hypotheses of Proposition~\ref{prop:main}. Hypothesis~\ref{hyp:lift} holds because its
lift is embedded (Corollary~\ref{cor:embedded}). Hypothesis~\ref{hyp:submersion} holds by Lemma~\ref{lem:reduction}
applied to Family $d$, since $\delta$ restricted to ${\cal P\/}$ is a
submersion at $\Psi(\zeta)$ by Lemma~\ref{subd}. So by Proposition~\ref{prop:main} there are embedded tori near $\Psi(\zeta)$
with $\Phi=(\zeta,0)$, that is, origami tori modeled on ${\cal T\/}$ with
modulus $\zeta$. As this holds for every $\zeta\in{\cal M\/}$, ${\cal T\/}$
is universal by \S\ref{sec:need}.

\begin{remark}
These same computations serve a more geometric proof, which needs no
smooth extension of the modulus. By Lemma~\ref{subd}, near our family the
flat tori form a $17$-dimensional manifold ${\cal F\/}$, and since
$\mu(\Psi(\xi))=\xi$, the modulus $\mu:{\cal F\/}\to\H^2$ is a submersion
there, with $15$-dimensional fibers. Vertical motions of a planar torus
change the edge lengths only to second order, so the lift $\Psi+tz$ is
tangent to the fiber through $\Psi(\zeta)$, and stays within $O(t^2)$ of
it, while the embeddedness determinants of \S\ref{sec:lift} grow linearly
in $t$. So a curve in the fiber leaving $\Psi(\zeta)$ with velocity $z$
differs from the lift by $O(t^2)$, and is embedded for small $t>0$: an
origami torus of modulus $\zeta$.
\end{remark}

\section{How We Found the Family}
\label{sec:found}

As described in \S\ref{sec:history}, our family grew out of numerical
experiments by F.L. and S.T., which searched all seven $8$-vertex
triangulations and found many approximate origami tori. F.L. took tori with
moduli very near the square and hexagonal lattices, drove them
numerically onto those lattices, and then down into the plane, keeping
them flat with fixed modulus. This did not always work, but when it did,
it produced flattened tori with simple incidences; on the hexagonal one,
for instance, the vertices $0,2,4,5,6$ lie on a circle
[{\bf L\/}, \S6.1]. Imposing these incidences as exact equations gave the
exact flattened tori of [{\bf L\/}].

F.L.'s incidences generalized to families of flattened tori along
vertical rays in ${\cal M\/}$. Trying to explain the incidences, we found
their source: they are what it takes for the folding to close up in the
Euclidean group of the plane. Solving these equations step by step, as in
\S\ref{sec:finding}, gives the full family of flattened tori with a given
folding pattern. The family of F.L.'s square example does not reach the
hexagonal lattice, but that of his hexagonal example covers ${\cal M\/}$.
It is the family of this paper.

Families that cover ${\cal M\/}$ seem to be rare. Up to symmetry, the seven $8$-vertex triangulations have $67{,}645$ folding
patterns, and three simple conditions rule out all but $2{,}940$. 
The signed angles of the triangles around each vertex must sum to $0$ or $\pm2\pi$, which constrains how many triangles can be face up.
No island may wrap around the
torus, since $\phi$ moves it by a single isometry. And every vertex must
lie on a crease, since the cone angle at a vertex inside an island cannot
change under horizontal motions, so $\delta$ restricted to ${\cal P\/}$
cannot be a submersion.
For the remaining $2{,}940$ patterns, we attempted to solve the folding equations numerically
at many moduli. In the end, only $18$ appeared to fold up across all of ${\cal M\/}$.

The existence of an embedded vertical lift is also quite restrictive: for
universality, even a single modulus where it fails rules a pattern out.
At a given flattened torus, the conditions describing embedded lifts cut out an open cone, and when this is nonempty it is usually easy to numerically exhibit a point.
So either a
pattern lifted quickly, and we moved on to other moduli, or it failed at
many hundreds of flattened tori of one modulus, and we discarded it.
Only two patterns survived: the pattern of F.L.'s hexagonal torus, on
${\cal T\/}(46666677)$, and one on ${\cal T\/}(55666677)$.

These two were our candidates for universality, and for both we then
confirmed everything exactly. To keep the polynomial checks simple, we
restrict to $2$-dimensional sections over ${\cal M\/}$, taking
$u,k,\lambda$ to be functions of $\zeta$, so that every polynomial is in
$x$ and $y$ rather than five variables. We chose the sections and the
heights together, by trial and error: we tried many candidate sections,
and for each searched for heights that work along the whole section. The
closer the sections and heights are to constant, the lower the degrees
of the polynomials in the final checks.

On ${\cal T\/}(46666677)$ the vertex coordinates are simple
(\S\ref{sec:need}), so the polynomials have low degree and every check
reduces to one in a single variable. On ${\cal T\/}(55666677)$ they are
much more complicated, and even the simplest sections we know lead to
polynomials that are irreducibly multivariable and of high degree: the
lift proof rests on the signs of $131$ irreducible polynomials in two
variables, of degree up to $10$ in one and $4$ in the other, and the
Jacobian has irreducible factors of degrees $(14,8)$ and $(22,11)$, the
latter with $220$ terms. So we wrote the paper around the simpler one.

\section{Polynomial Checks}
\label{sec:checks}

Covering ${\cal M\/}$ with families promoted the pointwise numerical checks
required by Proposition~\ref{prop:main} to familywise polynomial checks:
each says that some polynomial in $x$ and $y$ does not vanish on a strip
${\cal M\/}|_I=\{\zeta\in{\cal M\/}: x\in I\}$. Our two families were
chosen so that each check comes down to a polynomial of the form
$A(x)+B(x)\,y^2$ (see \S\ref{sec:found}), which is linear in $y^2$ along
each vertical line, and Lemma~\ref{lem:nonvanishing} below reduces it to
checks on polynomials in $x$ alone.  
We also describe a second method, strip dominance, which is the one used by
the file {\bf UniversalityProof\/}. It needs none of this structure: it
applies directly to the rational functions in our lemmas, without first
putting them in this form, and certifies each check exactly.

\begin{lemma}
\label{lem:nonvanishing}
Let $I\subset[0,\frac12]$ be an interval and $P(x,y)=A(x)+B(x)\,y^2$,
with $A$ and $B$ polynomials. Suppose that for every $x\in I$,
\begin{enumerate}
\item $A(x)+B(x)(1-x^2)\ne0$, and
\item $B(x)$ is either zero or has the same sign as $A(x)+B(x)(1-x^2)$.
\end{enumerate}
Then $P$ has no zero on ${\cal M\/}|_I$.
\end{lemma}

\startproof
Fix $x\in I$. In ${\cal M\/}|_I$ the vertical line through $x$ is the ray
$y^2\ge1-x^2$, and along it $P$ is linear in $y^2$, with slope $B(x)$.
At the bottom of the ray its value is $A(x)+B(x)(1-x^2)\ne0$, and the
slope is zero or has the same sign, so $P$ never vanishes on the ray.
\endproof

We now describe the general-purpose method. Its basic tool is a sign
test for polynomials on the region $\Omega=[0,1]\times[0,\infty)$.

\paragraph{Strip Dominance:}
Given a polynomial $g(t)=a_0+\cdots+a_nt^n$ we define
\begin{equation}
  \langle g \rangle = \min(A_0,\ldots,A_n), \hskip 20 pt
  A_k=a_0+\cdots+a_k.
\end{equation}
Then $g>0$ on $[0,1]$ if $\langle g\rangle>0$, and $g\ge0$ on $[0,1]$ if
$\langle g\rangle\ge0$, since
\begin{equation}
g(t)=(1-t)(A_0+A_1t+\cdots+A_{n-1}t^{n-1})+A_nt^n.
\end{equation}
We call $G(t,h)$ {\it strip dominant\/} if for one choice of sign
\begin{equation}
  \pm G(t,h)=g_0(t)+g_1(t)\,h+\cdots+g_m(t)\,h^m, \hskip 15 pt
  \langle g_0 \rangle>0, \hskip 15 pt
  \langle g_1\rangle,\ldots,\langle g_m \rangle \geq 0.
\end{equation}
A strip dominant $G$ does not vanish on $\Omega$. More generally, $G$
does not vanish on $\Omega$ if $G(t/2,h)$ and $G(1-t/2,h)$ are strip
dominant; we call this the tweaked criterion.

\paragraph{Rational Parametrization:}
To apply this to a strip ${\cal M\/}|_{[x_0,x_1]}$, which is bounded
below by an arc of the unit circle, we change coordinates to turn it
into $\Omega$. As $u$ goes from $0$ to $1$, the point
\begin{equation}
  \bigl(x(u),y(u)\bigr)=\Bigl(\frac{2u}{1+u^2},\ \frac{1-u^2}{1+u^2}\Bigr)
\end{equation}
runs along the unit circle from $(0,1)$ to $(1,0)$, with rational
coefficients. We choose rationals $r_0<r_1$ with
\begin{equation}
  \label{outer}
  x(r_0) \leq x_0<x_1 \leq x(r_1),
\end{equation}
so that the arc from $u=r_0$ to $u=r_1$ lies over $[x_0,x_1]$, and rescale
it to $t\in[0,1]$ by $u=(1-t)\,r_0+t\,r_1$. Then
\begin{equation}
  \label{rat0}
  F(t,h)=\bigl(x(u),\,y(u)+h\bigr)
\end{equation}
maps $\Omega$ onto the part of ${\cal M\/}$ above this arc, which
contains the strip. We use $r_0=0$ and $r_1=43/512$ for Family 1, which
covers $0\le x\le 0.16679$, and $r_0=143/1704$ and $r_1=144/512$ for
Family 2, which covers $0.1666666\le x\le 0.5213$. The choice of
$143/1704$ is extremely tight: just left of $x=\frac16$, at
$x\approx0.1666641$, the lift of Family 2 stops being embedded.

\paragraph{Certificate:}
Given a rational function $f$, we write $f\circ F=P/Q$ with $P$ and $Q$
integer polynomials in $t$ and $h$, and check that each of $P$ and $Q$
satisfies the strip dominance criterion or its tweaked version. Then
neither vanishes on $\Omega$, so $f$ is defined and nonzero on
${\cal M\/}|_{[x_0,x_1]}$.

\section{References}

\setlength{\parskip}{0.8\baselineskip}

\noindent
[{\bf ALM\/}] P. Arnoux, S. Leli\`evre, A. M\'alaga, {\it Diplotori:
  a family of polyhedral flat tori.\/} in preparation

\noindent
[{\bf Br\/}]  U. Brehm, Oberwolfach report (1978)

\noindent
[{\bf BZ1\/}] Y. D. Burago and V. A. Zalgaller,
{\it Polyhedral realizations of developments\/} (in Russian)
Vestnik Leningrad Univ. 15, pp 66--80 (1960)

\noindent
[{\bf BZ2\/}] Y. D. Burago and V. A. Zalgaller,
{\it Isometric piecewise-linear embeddings of
two-dimensional manifolds with a polyhedral
metric into $\R^3$\/},
Algebra i Analiz {\bf 7\/}(3), pp 76--95 (1995).
English translation: St. Petersburg Math. J.
{\bf 7\/}(3), pp 369--385 (1996).

\noindent
[{\bf DS\/}] P. G. Doyle and R. E. Schwartz,
{\it Collapsibility and Near Universality for Vertex Minimal Paper Tori\/},
arXiv:2510.12623 (2025)

\noindent
[{\bf L\/}] F. Lander,
{\it The Two Most Symmetric Flat Tori as Eight-Vertex Paper Tori\/},
arXiv:2609.10703 (2026)

\noindent
[{\bf LT\/}] F. Lazarus, F. Tallerie, {\it A Universal Triangulation
  for Flat Tori\/}, Discrete and Computational Geometry {\bf 71\/}
(2024) pp 278--307

\noindent
[{\bf S1\/}] R. E. Schwartz,  {\it Vertex Minimal Paper Tori\/},
arXiv 2507.14998  (2025)

\noindent
[{\bf S2\/}] R. E. Schwartz,  {\it The Most Efficient Origami
  Torus\/}, P.N.A.S. {\bf 123\/} (21) (2026)

\noindent
[{\bf S3\/}] R. E. Schwartz,
{\it The Projective Heat Map\/},
Mathematical Surveys and Monographs {\bf 219\/},
American Mathematical Society, 2017.

\noindent
[{\bf Se\/}] H. Segerman, {\it Visualizing Mathematics with 3D Printing\/},
Johns Hopkins U. Press  (2016)

\noindent
[{\bf T\/}] T. Tsuboi, {\it On origami embeddings of flat tori\/},
arXiv 2007.03434 (2020)

\noindent
[{\bf W\/}] Wolfram Research, Inc.,
{\it Mathematica\/}, Wolfram Research, Inc., Champaign, IL,
\texttt{https://www.wolfram.com/mathematica}

\end{document}